\documentstyle[titlepage,a4,12pt]{article}

\newtheorem{thm}{Theorem}
\newtheorem{lem}[thm]{Lemma}

\newcommand{\epr}{\hfill {$\Box$}}

\begin{document}
\begin{titlepage}
  \title {Towards a Characterisation of Pfaffian graphs}
  \author{Charles H. C. Little \footnote{The first author thanks the University of Klagenfurt for its hospitality while this research was undertaken.}\\ Massey University\\Palmerston North,
    New Zealand\\e-mail: c.little@massey.ac.nz \and Franz Rendl\\University of
    Klagenfurt\\Klagenfurt, Austria\\e-mail: franz.rendl@uni-klu.ac.at \and
Ilse Fischer\\University of Klagenfurt\\Klagenfurt, Austria\\e-mail: ilse.fischer@uni-klu.ac.at} \date\today
\end{titlepage}
\maketitle 
\begin{abstract}
A bipartite graph $G$ is known to be Pfaffian if and only if it does not contain an even subdivision $H$ of $K_{3, 3}$ such that $G-VH$ contains a $1$-factor. However a general characterisation of Pfaffian graphs in terms of forbidden subgraphs is currently not known. In this paper we describe a possible approach to the derivation of such a characterisation. We also extend the characterisation for bipartite graphs to a slightly more general class of graphs.
\end{abstract}
\newpage
\begin{section}
  {\Large\bf Introduction}
\end{section}
The graphs considered in this paper are finite and have no loops or
multiple edges. They are also undirected and connected unless an
indication to the contrary is given. If $v$ and $w$ are vertices in a
directed graph, then $(v, w)$ denotes an edge joining $v$ and $w$ and
directed from $v$ to $w$. If $G$ is any graph, then we denote its
vertex set by $VG$ and its edge set by $EG$. A {\it $1$-factor} of $G$
is a subset $f$ of $EG$ such that every vertex has a unique edge of
$f$ incident on it.

Let $G^*$ be a directed graph with an even number of vertices, and let
$F$ be the set $\{f_1, f_2,\cdots,f_k\}$ of $1$-factors of $G^*$. For
all $i$ write
\begin{center}
  $f_i=\{(u_{i1}, w_{i1}), (u_{i2}, w_{i2}),\ldots,(u_{in},
  w_{in})\}$,
\end{center}
where $u_{ij}, w_{ij} \in VG^*$ for all $j$. Associate with $f_i$ a
plus sign if
\begin{center}
  $u_{i1} w_{i1} u_{i2} w_{i2} \cdots u_{in} w_{in}$
\end{center}
is an even permutation of
\begin{center}
  $u_{11} w_{11} u_{12} w_{12} \cdots u_{1n} w_{1n}$,
\end{center}
and a minus sign otherwise. Thus the signs of the $1$-factors are
independent of the order in which their edges have been written. They
are dependent on the choice of $f_1$, but the resulting partition of
$F$ into two complementary subsets is not. If $G$ is an undirected
graph, we say that $G$ is a {\it Pfaffian} graph if there exists a
directed graph $G^*$ with vertex set $VG$ and edge set $EG$ such that
all the $1$-factors of $G^*$ have the same sign. We say that $G^*$ is
a {\it Pfaffian} orientation for $G$.

It is a tantalising problem to characterise Pfaffian graphs in terms
of forbidden subgraphs. Pfaffian bipartite graphs have been so
characterised by Little \cite{Li:75}. As a further contribution to a general
characterisation, the present paper represents an attempt to
strengthen this result.

Further progress may be possible because of a theorem of Lov\'{a}sz
and Plummer \cite[Theorem 5.4.6]{LoPl:86} on ear decompositions of $1$-extendible graphs. In
order to describe this theorem, we need some more definitions.

A graph is {\it $1$-extendible} if every edge has a $1$-factor
containing it. We identify paths and circuits with their edge sets. If
$X$ is a path or circuit, then we denote by $VX$ the set of vertices
of $X$. If $P$ is a path and $u$, $v \in VP$, then we denote by $P [u,
v]$ the subpath of $P$ joining $u$ to $v$. An {\it ear} is a path of
odd cardinality. A circuit is {\it alternating} with respect to two
given $1$-factors if it is included in their symmetric difference. A
circuit that is alternating with respect to a $1$-factor $f$ is also
said to be {\it $f$-alternating}, or {\it consanguineous} (with
respect to $f$). A path $P$ is {\it $f$-alternating} if every internal
vertex of $P$ is incident with an edge of $P \cap f$.

Now let $H$ be a $1$-extendible subgraph of a $1$-extendible graph
$G$. Let $A$ be an alternating circuit in $G$ which includes $EG$-$EH$
and meets $EH$. Then an $AH$-{\it arc} (or an $HA$-{\it arc}) is a subpath of $A \cap EH$ of
maximal length, and an $A \bar H$-{\it arc} (or an $\bar HA$-{\it arc}) is a
subpath of $EG$-$EH$, of maximal length, whose internal vertices are
in $VG$-$VH$. If there are $n$ such $A \bar H$-arcs, and each is an
ear, then we say that $G$ is obtained from $H$ by an {\it $n$-ear
  adjunction}. An {\it ear decomposition} of $G$ is a sequence $G_0,
G_1,\cdots,G_t$ of $1$-extendible subgraphs of $G$ such that $G_0$ is
isomorphic to $K_2$, $G_t=G$ and, for each $i>0$, $G_i$ is obtained
from $G_{i-1}$ by an $n$-ear adjunction with $n=1$ or $n=2$. The
theorem of Lov\'{a}sz and Plummer alluded to earlier asserts that
every $1$-extendible graph has an ear decomposition. It can be stated
as follows.  \newline  \newline {\bf Theorem 1} {\it Let $f$ be a
  $1$-factor in a $1$-extendible graph $G$. Let $H$ be a
  $1$-extendible proper subgraph of $G$ such that $EH \not= \emptyset$
  and $f \cap E H$ is a $1$-factor of $H$. Then $G$ contains an
  $f$-alternating circuit $A$ that admits just one or two $A \bar
  H$-arcs.}  \newline \newline In the special case where $G$ is
bipartite, it had already been shown that $A$ may be chosen to
admit just one $A \bar H$-arc.  In fact this result was implicit in \cite{Li:75}, and was used there to prove the following characterisation of
Pfaffian bipartite graphs.  \newline \newline  {\bf Theorem 2} {\it A
  bipartite graph $G$ is non-Pfaffian if and only if $G$ contains an
  even subdivision $J$ of $K_{3,3}$ such that $G-VJ$ has a
  $1$-factor.}  \newline  \newline Here we need to explain the term
``even subdivision''. An {\it edge subdivision} of a graph $G$ is
defined as a graph obtained from $G$ by replacing an edge joining
vertices $v$ and $w$ with a path $P$ joining $v$ and $w$ such that $VP
\cap VG = \{v, w\}$. The edge subdivision is {\it even} if $|VP|$ is
even. A graph $H$ is a {\it subdivision} of $G$ if for some positive
integer $k$ there exist graphs $G_0, G_1,\cdots,G_k$ such that
$G_0=G$, $G_k=H$ and, for all $i>0$, $G_i$ is an edge subdivision of
$G_{i-1}$. If $G_1, G_2,\cdots,G_k$ can be chosen so that in addition
$G_i$ is an even edge subdivision of $G_{i-1}$ for all $i>0$, then $H$
is said to be an {\it even subdivision} of $G$.

The idea behind the proof of Theorem $2$ runs as follows. Clearly we
may assume that $G$ is $1$-extendible. If $G$ contains $J$, then $G$
is easily seen to be non-Pfaffian. Suppose on the other hand that $G$
is non-Pfaffian. We construct an ear decomposition $G_0,
G_1,\cdots,G_t$ of $G$. Since $G$ is bipartite, we may assume that,
for each $i>0$, $G_i$ is obtained from $G_{i-1}$ by the adjunction of
a single ear. As $G_0$ is Pfaffian but $G$ is not, there exists a
smallest positive integer $j$ such that $G_j$ is non-Pfaffian. The
graph $G_j$ is studied in detail and eventually shown to contain $J$.

Theorem $1$ provides a possible way to generalise this argument. If we
drop the assumption that $G$ is bipartite, then $G_j$ is obtained from
$G_{j-1}$ by the adjunction of one or two ears. Morever $G_{j-1}$
might or might not be bipartite. We therefore distinguish four cases:
\begin{enumerate}
\item[(a)] a $1$-ear adjunction to a bipartite graph;
\item[(b)] a $2$-ear adjunction to a bipartite graph;
\item[(c)] a $1$-ear adjunction to a non-bipartite graph;
\item[(d)] a $2$-ear adjunction to a non-bipartite graph.
\end{enumerate}
In each case the goal is to find a member of a class of forbidden
graphs as a subgraph of $G_j$. In case (a), where $G_j$ is obtained
from the bipartite graph $G_{j-1}$ by the adjunction of a single ear,
this goal was achieved in \cite{Li:75}, by showing that $G_j$ contains an even
subdivision of $K_{3,3}$. The remaining cases are still open.

It turns out that it is unnecessary to investigate case (d). In order
to explain the reason, we need first to provide an account of the
brick decomposition procedure of Edmonds, Lov\'{a}sz and Pulleyblank \cite{EdLoPu:82}. A graph $G$ is {\it bicritical} if $G-\{u, v\}$ has a $1$-factor
for any pair of vertices $u$ and $v$. A {\it brick} is a $3$-connected
bicritical graph.

Suppose that $G$ is a graph which is $1$-extendible but not
bicritical. Then (see \cite[Theorem 5.2.2(d)]{LoPl:86}) it has a maximal set $S$ of vertices such
that $|S|\ge 2$ and $G-S$ has exactly $|S|$ odd components (components
with an odd number of vertices). Let $|S|=k$, and let $H_1,
H_2,\cdots,H_k$ be the odd components of $G-S$. For each $i$ let $G_i$
be the graph obtained from $G$ by contraction of the subgraph
$G-VH_i$. Let $G_0$ be the bipartite graph obtained from $G$ by the
successive contraction of $H_i$ for each $i$. We call $G_0$ the {\it
  frame}. It is shown in \cite[Theorem 5.2.6]{LoPl:86} that $G_0, G_1,\cdots,G_k$ are
$1$-extendible. It is also clear that, for any $i$, the graph obtained
from $G$ by the contraction of $H_i$ is $1$-extendible.

Now we discard $G_0$ and any of $G_1, G_2,\cdots,G_k$ that are
isomorphic to $K_2$, and we file those of $G_1, G_2,\cdots,G_k$ that
are bicritical. The procedure is then repeated recursively for the
remaining graphs. Eventually a family of bicritical graphs is
obtained.

A bicritical graph may be decomposed into bricks as follows. Let $G$
be bicritical but not a brick. It follows that the connectivity of $G$
is $2$. Hence there are vertices $u, v$ such that $G-\{u,v\}$ is not
connected. Let $G_1', G_2',\cdots,G_l'$ be the components of $G-\{u,
v\}$, and for each $i$ let $G_i$ be the graph obtained from $G[VG_i'
\cup \{u, v\}]$ by adjoining an edge between $u$ and $v$ if they are
not already adjacent. It is shown in \cite[Lemma 5.2.8]{LoPl:86} that $G_i$ is bicritical.
This procedure is repeated recursively until a list of bricks is
obtained.

The brick decomposition procedure described above motivates the study
of two operations. Firstly, let $v$ and $w$ be vertices, of equal
degree $d$, in graphs $H$ and $K$ respectively. Let $v_1,
v_2,\cdots,v_d$ be the neighbours of $v$ and $w_1, w_2,\cdots,w_d$
those of $w$. Let $G$ be the graph obtained from $(H-\{v\}) \cup
(K-\{w\})$ by adjoining an edge between $v_i$ and $w_i$ for each $i$.
Then we say that $G$ is formed by {\it splicing} $H$ and $K$ at $v$
and $w$ respectively. We call $EG-(EH \cup EK)$ the {\it splice} of
$G$.

Secondly, let $x$ and $y$ be edges in graphs $H$ and $K$ respectively.
Let $L$ be the graph obtained from $H$ and $K$ by identifying $x$ and
$y$ to form an edge $e$. Then graphs $L$ and $L-\{e\}$ are said to be
obtained from $H$ and $K$ by {\it gluing} $H$ and $K$ at $x$ and $y$.
The brick decomposition procedure shows that any $1$-extendible graph
may be constructed from bricks by gluing and splicing, where the
graphs being glued and spliced each have more than one edge.

The effect of the operations of gluing and splicing on the Pfaffian
property has been studied in \cite{LiRe:91}, and the results are as follows.
Splicing or gluing a non-Pfaffian graph to another graph yields a
non-Pfaffian graph. Gluing Pfaffian graphs yields a Pfaffian graph. In
regard to the splicing of $1$-extendible Pfaffian graphs, however, the
situation is a little more complicated. Let us define a {\it cut} in a
graph $G$ to be a minimal set $X$ of edges such that $G-X$ has more
components than has $G$. An example is a splice of $G$. A cut $X$ is
called {\it tight} if $| f \cap X| = 1$ for every $1$-factor $f$. It
is shown in \cite{LiRe:91} that if $G$ is obtained by splicing $1$-extendible
Pfaffian graphs, and the splice in $G$ is tight, then $G$ is Pfaffian.
It is easily checked that, if the brick decomposition procedure is
reversed, then all splices are tight.

From the discussion in the previous paragraph, we see that in order to
characterise Pfaffian graphs it suffices to characterise Pfaffian
bricks.

Now we introduce a theorem, due to Carvalho, Lucchesi and Murty \cite{CaLuMu:98},
concerning the smallest number of $2$-ear adjunctions in an ear
decomposition. The {\it cycle space}, ${\cal C}(G)$, of a graph $G$ is
the vector space spanned by the circuits of $G$, where the sum of
vectors is defined as their symmetric difference. The {\it alternating
  space}, ${\cal A}(G)$, of $G$ is the subspace of ${\cal C}(G)$
spanned by the alternating circuits. Note that the total number of
ears adjoined in the course of an ear decomposition of a
$1$-extendible graph $G$ is $\dim {\cal C}(G)$, for if $G$ is obtained
from a subgraph $H$ by an $n$-ear adjunction then
\begin{center}
  $\dim {\cal C}(G) - \dim {\cal C}(H) = n$.
\end{center}
On the other hand,
\begin{center}
  $\dim {\cal A}(G) - \dim {\cal A}(H) \ge 1$.
\end{center}
These results imply that a lower bound for the number of $2$-ear
adjunctions in an ear decomposition is $\dim {\cal C}(G) - \dim {\cal
  A}(G)$. The theorem of Carvalho, Lucchesi and Murty alluded to
earlier is that this lower bound can always be met.  \newline \newline
{\bf Theorem 3} {\it The minimum number of $2$-ear adjunctions in an
  ear decomposition of a $1$-extendible graph $G$ is $\dim {\cal C}(G)
  - \dim {\cal A}(G)$.}  \newline \newline In the case of a bipartite
graph, it has been shown that no $2$-ear adjunctions are
necessary. It follows that if $G$ is $1$-extendible and bipartite,
then ${\cal A}(G) = {\cal C}(G)$. In fact, if $G_0, G_1,\cdots,G_t$ is
an ear decomposition of a graph $G$, then $G$ is bipartite if and only
if, for each $i>0$, $G_i$ is obtained from $G_{i-1}$ by the adjunction
of a single ear.

The value of $\dim {\cal C}(G) - \dim {\cal A}(G)$ has been given by
Lov\'{a}sz \cite{Lo:87}. In the case where $G$ is a brick, this number is $2$
if $G$ is isomorphic to the Petersen graph and $1$ if $G$ is any other
brick. After noting that the Petersen graph is non-Pfaffian, we can
therefore focus our attention on $1$-extendible graphs $G$ for which $$\dim {\cal
  C}(G) - \dim {\cal A}(G) = 1.$$
By Theorem 3, there is an ear
decomposition $G_0, G_1,\cdots,G_t$ of $G$ for which there is a unique
$j>0$ such that $G_j$ is obtained from $G_{j-1}$ by a $2$-ear
adjunction. Then $G_{j-1}$ is bipartite but $G_j$ is not, since $j$ is
the smallest positive integer such that $G_j$ is obtained from
$G_{j-1}$ by a $2$-ear adjunction. Hence $G_j$ satisfies case (b). If
$0<i<j$, then $G_i$ is bipartite and therefore satisfies case (a). If
$j<i \le t$, then $G_i$ is non-bipartite and hence satisfies case (c).
Thus case (d) does not need to be considered.

The present paper addresses case (b), which arises as follows. Suppose that $G$ is a non-Pfaffian $1$-extendible graph
with ear decomposition $G_0, G_1,\cdots, G_t$, and that $t$ is the
smallest integer $j$ such that $G_j$ is non-Pfaffian, and also the
smallest integer $j$ such that $G_j$ is non-bipartite. Kasteleyn \cite{Ka:67}
has shown that the $1$-factors of a directed graph all have equal sign
if and only if all the alternating circuits are clockwise odd. (The
clockwise parity of a circuit of even length is the parity of the
number of its edges that are directed in agreement with a specified
sense.) Since $G_{t-1}$ is Pfaffian, we may assume it to be oriented
so that each of its alternating circuits is clockwise odd. Now extend
this orientation to $G_t$. Since $G_t$ is non-Pfaffian, there must be
alternating circuits $A$ and $B$ in $G_t$ of opposite clockwise
parity. As $G_{t-1}$ is bipartite but $G_t$ is not, $G_t$ is obtained
from $G_{t-1}$ by a $2$-ear adjunction. We show that $A$ and $B$ may
be chosen so that each traverses both ears and there are just one or
two $\bar AB$-arcs. We then proceed to consider the former case. A
sufficient condition for this case to arise is that $G [C \cup D]$ be
bipartite whenever $C$ and $D$ are alternating circuits that traverse
both ears and have the property that there are just two $\bar
CD$-arcs. Hence in this paper we confine our attention to graphs that
satisfy this condition.

It is shown in \cite{Li:75} that $K_{3,3}$ is non-Pfaffian. It follows that no
even subdivision of $K_{3,3}$ is Pfaffian. Our first lemma shows that
a graph $G$ is non-Pfaffian if it has a circuit $X$, of odd length,
such that the graph obtained from $G$ by contracting $VX$ is an even
subdivision of $K_{3,3}$. In general, let us say that a graph $G$ is
{\it reducible} to a graph $H$ if $G$ has a circuit $X$, of odd
length, such that $H$ is obtained from $G$ by contracting $VX$. Thus
any graph that is reducible to an even subdivision of $K_{3,3}$ is
non-Pfaffian. In fact, a graph $G$ must be non-Pfaffian if it has a
subgraph $H$ that is reducible to an even subdivision of $K_{3,3}$ and
has the property that $G-VH$ has a $1$-factor.

Now let $G$ be a $1$-extendible graph that can be obtained from a
$1$-extendible bipartite graph by a $2$-ear adjunction and suppose
that $G[C \cup D]$ is bipartite whenever $C$ and $D$ are alternating
circuits that traverse both ears and have the property that there are
just two $\bar CD$-arcs. In this paper we shall show that $G$ is
non-Pfaffian if and only if it has a subgraph $H$, reducible to an
even subdivision of $K_{3,3}$, such that $G-VH$ has a $1$-factor.

\begin{section}
  {\Large\bf Preliminary Lemmas}
\end{section}

A set $S$ of alternating circuits in a directed graph $G$ is called
{\it intractable} if each edge of $G$ belongs to an even number of
alternating circuits in $S$ and an odd number of the members of $S$
are clockwise even. The former property implies that the latter is
independent of the orientation of $G$. It is shown in \cite{Li:73} that $G$ is
Pfaffian if and only if it has no intractable set of alternating
circuits. In fact it suffices to consider only sets of alternating
circuits that are consanguineous with respect to a fixed $1$-factor.
\newline
\begin{lem}
  Let $G$ be a graph with a circuit $X$ of odd length. Let $H$ be the
  graph obtained from $G$ by contracting $VX$. If $H$ is not Pfaffian,
  then neither is $G$.
\end{lem}

{\bf Proof.} Let $v$ be the vertex $VX$ in $H$, and let $e_0,
e_1,\cdots,e_l$ be the edges incident on it. Since $H$ is
non-Pfaffian, it has a $1$-factor $f$ such that there is an
intractable set $S$ of $f$-alternating circuits. We must construct an
intractable set in $G$.

Without loss of generality, let $e_0 \in f$. Let $A$ be an
$f$-alternating circuit in $S$, and suppose that $A$ contains edges
$e_0$ and $e_i$ for some $i>0$. In $G$ the edges $e_0$ and $e_i$ are
incident on two vertices, $x$ and $y$ respectively, of $X$. Let $a$ be
the edge of $X$ that is incident on $y$ but not in $g$, where $g$ is
the unique $1$-factor of $G$ such that $g-f \subset X$. Let
$X'=X-\{a\}$, and let $A'=A \cup X' [x, y]$. (In other words, $A'$ is the path in $A$ of even length joining $x$ and $y$.) On the other hand, if
$e_0 \notin A$ then define $A'=A$. 

Let $S = \{S_1, S_2,\cdots,S_n\}$. Then the required intractable set
of $g$-alternating circuits in $G$ is $\{S'_1, S'_2,\cdots,S'_n\}$ where, for each $i$, $S'_i$ is constructed from $S_i$ as described above.
\epr \newline

We shall also need the following two lemmas, which give properties
concerning the structure of the sum and union, respectively, of two
consanguineous alternating circuits.  \newline
\begin{lem}
  Let $A_1, A_2$ be $f$-alternating circuits in a directed graph $G$
  with $1$-factor $f$. Then $A_1$ and $A_2$ are of opposite clockwise
  parity if and only if $A_1 + A_2$ includes an odd number of
  clockwise even alternating circuits.
\end{lem}

{\bf Proof.}  Let $A_1 = f + f_1$ and $A_2 = f + f_2$ for some
1-factors $f_1$ and $f_2$. Suppose that $A_1$ and $A_2$ are of
opposite clockwise parity. Without loss of generality, let $A_1$ be
clockwise odd and $A_2$ be clockwise even. Then $f$ has the same sign
as $f_1$ but the opposite sign from $f_2$. Hence $f_1$ and $f_2$ have
opposite sign. Since $A_1+A_2 = f_1 + f_2$, it follows that $A_1 +
A_2$ includes an odd number of clockwise even alternating circuits.
Similarly if $A_1$ and $A_2$ have the same clockwise parity, then
$f_1$ and $f_2$ have equal sign and $A_1 + A_2$ includes an even
number of clockwise even alternating circuits. \epr \newline
\begin{lem}
  Let $f$ be a 1-factor in a 1-extendible directed graph $G$.  Let $A$
  and $B$ be $f$-alternating circuits in $G$, of opposite clockwise
  parity, containing distinct independent edges $e_1$ and $e_2$ such
  that $e_1 \notin f$ and $e_2 \notin f$. Suppose that $G - \{ e_1\} $
  and $G - \{e_2\}$ are not bipartite but that $G-\{e_1,e_2\}$ is.
  Then $A \cup B$ includes alternating circuits $X$ and $Y$, of
  opposite clockwise parity and consanguineous with respect to some
  1-factor that contains neither $e_1$ nor $e_2$, such that there are
  just one or two $XY$-arcs, each $XY$-arc contains $e_1$ or $e_2$ and
  their union contains both.
\end{lem}

{\bf Proof.}  Let $G_0=G[A \cup B]$.  Since $A$ and $B$ are
$f$-alternating, they are $f_0$-alternating in $G_0$, where $f_0$ is
the 1-factor $f \cap (A \cup B)$ in $G_0 - \{e_1, e_2\}$. Thus $G_0$
is 1-extendible.  Let $A_0=A$ and $B_0=B$.

For some $i \geq 0$, let $G_i = G[ A_i \cup B_i]$ where, for some
1-factor $f_i$ of $G_i - \{e_1, e_2\}$, the circuits $A_i$ and $B_i$
are $f_i$-alternating circuits in $G_i$, of opposite clockwise parity,
containing $e_1$ and $e_2$.  Thus $G_i$ is 1-extendible. By the 2-ear
theorem applied to the subgraph $G [A_i]$, the graph $G_i$ contains an $f_i$-alternating circuit $B_i'$ such that
there are just one or two $\bar{A_i}B_i'$-arcs, and in the latter case
there is no $f_i$-alternating circuit $B_i^*$ such that there is only
one $\bar{A_i}B_i^*$-arc.

\medskip {\bf Case 1:} Suppose first that there is a unique
$\bar{A_i}B_i'$-arc $P_i$.  Let $P_i$ join vertices $u_i$ and $v_i$.
Since $A_i$ and $B_i'$ are $f_i$-alternating, $A_i + B_i'$ is an
alternating circuit $C_i$ and hence of even length. The hypotheses on
$G$ show that $C_i$ contains both $e_1$ and $e_2$ or neither, as any
circuit containing only one of $e_1$ and $e_2$ must be of odd length.

\medskip {\bf Subcase 1.1:} Suppose that $C_i$ contains $e_1$ and
$e_2$. Then $e_1 \notin B'_i$ and $e_2 \notin
B'_i$. Let $A_{i+1}= C_i$ and $B_{i+1}= B_i + B_i'$, so that $A_{i+1}$ and $B_{i+1}$ contain $e_1$ and $e_2$. Moreover $A_{i+1}+B_{i+1} = C_i + B_i + B'_i = A_i + B_i$, and so $A_{i+1}$
and $B_{i+1}$ are of opposite clockwise parity. Furthermore, since $A_{i+1}=A_i+B'_i = f_i+g_i+B'_i$ and $B_{i+1} = B_i + B'_i = f_i + h_i + B'_i$ for some
$1$-factors $g_i$ and $h_i$ of $G$, we see that $A_{i+1}$ and $B_{i+1}$ are
alternating with respect to the $1$-factor $f_{i+1}=f_i + B_i'$, which
contains neither $e_1$ nor $e_2$. Note also that the edge of $A_i \cap
B_i'$ incident on $u_i$ belongs to $f_i$ and therefore to neither
$A_{i+1}$ nor $B_{i+1}$. Thus if we define $G_{i+1}=G[A_{i+1} \cup
B_{i+1}]$ then $|EG_{i+1}| < |EG_i|$.

\medskip {\bf Subcase 1.2:} Suppose that $C_i$ contains neither $e_1$
nor $e_2$. If $C_i$ is clockwise even, then $A_i$ and $B_i'$ are the
required circuits $X$ and $Y$, by Lemma 2. Suppose therefore that
$C_i$ is clockwise odd.  Let $A_{i+1}=B'_i = A_i + C_i$, $B_{i+1}=B_i$
and $f_{i+1}=f_i$.  Then $A_{i+1}$ and $B_{i+1}$ are
$f_{i+1}$-alternating and contain $e_1$ and $e_2$. Since $C_i = A_i + A_{i+1}$ and $C_i$ is clockwise odd, it follows from Lemma 2 that $A_i$ and $A_{i+1}$ have the same clockwise parity. Hence $A_{i+1}$ and
$B_{i+1}$ are of opposite clockwise parity.  Moreover the edge of $C_i
- P_i$ incident on $u_i$ belongs to neither $A_{i+1}$ nor $B_{i+1}$,
and so $|EG_{i+1}|<|EG_i|$, where $G_{i+1}=G[A_{i+1} \cup B_{i+1}]$.

\medskip {\bf Case 2:} Suppose there are two $\bar{A_i}B_i'$-arcs,
$Q_i$ and $R_i$. Let $q_1$ and $q_2$ be the ends of $Q_i$, and let
$r_1$ and $r_2$ be the ends of $R_i$. If there exists a path $S$,
included in $A_i$, joining $q_1$ and $q_2$ and passing through both
$r_1$ and $r_2$, then $S \cup Q_i$ is an $f_i$-alternating circuit
which includes $Q_i$ but not $R_i$. This contradiction shows that each
path included in $A_i$ and joining $q_1$ and $q_2$ must pass through
exactly one of $r_1$ and $r_2$.  We may therefore assume without loss
of generality that
$$B_i' = M[r_1,q_1] \cup Q_i \cup N[q_2, r_2] \cup R_i,$$
where $M$ is
the path included in $A_i$ which joins $q_1$ to $q_2$ and passes
through $r_1$, and $N$ is the path $A_i - M$. Since $A_i$ and $B_i$
are $f_i$-alternating, the circuits $Q_i \cup M$, $Q_i \cup N$, $R_i \cup M [r_1, q_1] \cup N [q_1, r_2]$ and $R_i \cup M [r_1, q_2] \cup N [q_2, r_2]$ are
of odd length. Thus we may assume without loss of generality that either $e_1
\in M[ r_1, q_1]$ and $e_2 \in N[q_2, r_2]$ or $\{e_1, e_2\} \cap B'_i = \emptyset$.

Once again, let $C_i$ be the alternating circuit $A_i + B_i'$. If $C_i$ contains $e_1$ and $e_2$, then the argument proceeds as in Subcase 1.1. In the remaining case we have $\{e_1, e_2\} \subseteq B'_i$, and the argument proceeds as in Subcase 1.2. 

By the finiteness of $G$, there exists $j$ such that $G_{j+1}$ is not
defined. Then $A_j$ and $B_j'$ are the required circuits.  \epr
\newline \newline In the remainder of this section we let $e_1$ and
$e_2$ be distinct independent edges in a graph $G$ such that neither
$G-\{e_1\}$ nor $G-\{e_2\}$ is bipartite but $G-\{e_1, e_2\}$ is
bipartite and $1$-extendible. Let $f$ be a $1$-factor of $G-\{e_1,
e_2\}$. We suppose that there exists an $f$-alternating circuit $A$ in
$G$ containing $e_1$ and $e_2$, and we let $A'=A-\{e_1\}$. Let $e_1$
join vertices $u_1$ and $v_1$, and let $e_2$ join vertices $u_2$ and
$v_2$, where $e_2 \in A' [u_1, v_2]$.  \newline
\begin{lem}
  Let $C$ be an $f$-alternating circuit in $G-\{e_1, e_2\}$, and let a
  be an edge of $A' [v_1, v_2] \cap C$ joining vertices $u$ and $v$,
  where $a \in A' [v_2, u]$. Let $C' = C-\{a\}$. Let $b$ be an edge of
  $C'$ joining vertices $w$ and $x$, where $b \in C' [v, x]$.
   \begin{flushleft}
     (a) If $ b \in A'[v_1, v_2]$, then $b \in A'[v_2, w]$. \\
     (b) If $b \in A'[u_1, u_2]$, then $b \in A'[u_2, w]$.
   \end{flushleft}
\end{lem}

{\bf Proof.} Let $X$ be an $\bar A C$-arc joining vertices $y$ and
$z$. Since $A$ and $C$ are both $f$-alternating, there must be an edge
$c \in f \cap A \cap C$ incident on $y$ and an edge $d \in f \cap A
\cap C$ incident on $z$. If $e_2 \in A' [y, z]$ then the circuit $X
\cup A'[y, z]$ contains $e_2$ but not $e_1$, and is therefore of odd
length by the hypotheses on $G$. In this case exactly one of $c$ and
$d$ is in $A'[y, z]$. Similarly if $e_2 \not\in A'[y, z]$ then either
$c$ and $d$ are both in $A'[y, z]$ or neither of them is. The lemma
follows by applying these observations sequentially to all the $\bar A
C$-arcs.  \epr \newline

\begin{section}
{\Large\bf The Main Results}
\end{section}
In order to prepare for the proof of our first theorem, we need the
following rather technical definition.

Let $G$ be a graph with $1$-factor $f$. Let $A$ be an $f$-alternating
path in $G$ joining vertices $u$ and $v$. Let $x_0$ and $y_0$ be distinct vertices of $VA$, where $x_0 \in VA [y_0, v]$. Suppose that the edges of $f$ incident on $x_0$ and $y_0$ belong to $A[x_0, v]$ and $A[u, y_0]$ respectively. Let $C_1, C_2,\cdots,C_k$ be
$f$-alternating circuits with the following properties.  
\newline

\begin{enumerate} 
\item[(a)] For each $i$ there exists a unique $AC_i$-arc. Let $x_i$ and $y_i$ be the ends of this arc, where $x_i \in VA [u, y_i]$.
\item[(b)] We have $x_1 \in VA [u, y_0]-\{y_0\}$, $y_1 \in VA [y_0, v]-\{y_0\}$ and
  $y_k \in VA [x_0, v]-\{x_0\}$.
\item[(c)] If $k > 1$, then for each $i \ge 2$ we have $x_i \in VA
  [y_{i-2}, y_{i-1}]-\{y_{i-2}, y_{i-1}\}$ and $y_i \in VA [y_{i-1}, v]-\{y_{i-1}\}$.
\end{enumerate}
Then the sequence $$(C_1, C_2,\cdots,C_k)$$
is called a {\it cascade} of {\it length} $k$ along $A [y_0, x_0]$ from $y_0$ to $x_0$.
\newline 
\newline 
{\bf Theorem 4} {\it Let $G$ be a
  $1$-extendible graph with $1$-factor $f$. Let $e_1$ and $e_2$ be
  distinct independent edges of $EG-f$ such that neither $G-\{e_1\}$
  nor $G-\{e_2\}$ is bipartite but $G-\{e_1, e_2\}$ is bipartite, Pfaffian and
  $1$-extendible. Suppose there exist $f$-alternating circuits $A$ and
  $B$, both containing $e_1$ and $e_2$, such that there is a unique
  $AB$-arc and $A+B$ is clockwise even under a Pfaffian orientation of $G-\{e_1, e_2\}$. Then $G$ has a subgraph $H$,
  reducible to an even subdivision of $K_{3,3}$, such that $G-VH$ has
  a $1$-factor.}  \newline

{\bf Proof.} Let $e_1$ join vertices $u_1$ and $v_1$, and let
$A'=A-\{e_1\}$. Let $e_2$ join vertices $u_2$ and $v_2$, where $u_2
\in VA' [u_1, v_2]$. Let $x_0$ and $y_0$ be the ends of the unique
$\bar AB$-arc $P$. Since $\{e_1, e_2\} \subset A \cap B$, we may
assume with no less generality that $x_0 \in VA' [v_1, v_2]$ and $y_0
\in VA' [v_2, x_0]$. Moreover the edge $a$ of $f$ incident on $y_0$
is in $A' [v_2, y_0]$ and the edge of $f$ incident on $x_0$ is in $A'
[x_0, v_1]$. We may assume $A$ and $B$ to have been chosen to maximise
$|A' [v_2, y_0]| + |A'[x_0, v_1]|$. Let $b$ be the edge of $A' [x_0,
y_0]$ incident on $y_0$.

\medskip {\bf Claim 1:} There is a cascade ${\cal C}$ along $A' [y_0, x_0]$ from $y_0$ to $x_0$ such that distinct elements of this cascade meet, if at all, only in $A' [x_0, y_0]$.

Let $K = G - \{e_1, e_2\}$. Since $K$ is $1$-extendible, there is an $f$-alternating circuit $C_1$ in $K$ containing $a$ and $b$. Let $C'_1 = C_1 - \{a\}$. 

Since $A' [v_2, y_0] \cap C_1 \not= \emptyset$ and $A' [y_0, v_1] \cap C_1 \not= \emptyset$, we may select a vertex $$x_1 \in (VA' [v_2, y_0] - \{y_0\}) \cap VC_1$$ that minimises $|C'_1 [x_1, y_0]|$ and a vertex $$y_1 \in VA' [y_0, v_1] \cap VC_1$$ that minimises $|C'_1 [x_1, y_1]|$. We may assume that $C_1$ has been chosen to minimise $|A' [y_1, v_1]|$, and by Lemma 4 we may assume that $A' [x_1, y_1]$ is an $AC_1$-arc.

From Lemma 4 we see that $y_1 \not= x_0$. If $y_1 \in VA' [y_0, x_0]$, then we repeat the argument. Thus there must be an $f$-alternating circuit $C_2$ in $K$ containing both the edges of $A$ incident on $y_1$. Let $C'_2 = C_2 - \{a'\}$, where $a'$ is the edge of $f$ incident on $y_1$. We may select a vertex $$x_2 \in (VA' [v_2, y_1] - \{y_1\}) \cap VC_2$$ to minimise $|C'_2 [x_2, y_1]|$ and a vertex $$y_2 \in VA' [y_1, v_1] \cap VC_2$$ to minimise $|C'_2 [x_2, y_2]|$. In fact it follows from the choice of $C_1$ that $x_2 \in VA' [y_0, y_1]$. We may assume that $C_2$ has been chosen to minimise $|A' [y_2, v_1]|$. By Lemma 4 we may also suppose that $A' [x_2, y_2]$ is an $AC_2$-arc.

Suppose that $$C'_1 [x_1, y_1] \cap C'_2 [x_2, y_2] \not= \emptyset.$$ Let $w$ be the vertex of $VC'_1 [x_1, y_1] \cap VC'_2 [x_2, y_2]$ that minimises $|C'_2 [w, y_2]|$. Since the circuit $$C'_1 [w, y_1] \cup A' [y_1, y_2] \cup C'_2 [y_2, w]$$ contains neither $e_1$ nor $e_2$ and is therefore of even length, it follows that $$C'_1 [x_1, w] \cup C'_2 [w, y_2] \cup A' [y_2, x_1]$$ is an $f$-alternating circuit containing $a$ and $b$. Since $y_2 \in VA' [y_1, v_1] - \{y_1\}$, the choice of $C_1$ is contradicted. We conclude that $C'_1 [x_1, y_1] \cap C'_2 [x_2, y_2] = \emptyset$. 

From Lemma 4 we see that $y_2 \not= x_0$. If $y_2 \in VA' [y_0, x_0]$, then we repeat the argument inductively. By the finiteness of $G$, there exists $k > 0$ such that $y_k \in VA' [x_0, v_1]$. Then $(C_1, C_2,\cdots, C_k)$ is the required cascade ${\cal C}$ along $A' [y_0, x_0]$ from $y_0$ to $x_0$. The proof of the claim is now complete. \medskip

For each $i$, let $p_i$ be the number of $BC_i$-arcs included in $P$, and let $q_i$ be the number of $AC_i$-arcs. We may assume $f$, $A$, $B$ and ${\cal C}$ chosen to minimise 
\begin{eqnarray}
\sum_{i=1}^k (p_i + q_i).
\end{eqnarray}

Since distinct elements of ${\cal C}$ meet, if at all, only in $A' [x_0, y_0]$, it follows that $$\bigcup_{i=1}^{k-1} (C'_i [x_i, y_i] \cup A' [y_i, x_{i+1}]) \cup C'_k [x_k, y_k]$$ is a path $X$ joining $x_1$ to $y_k$. Traversed from $x_1$ to $y_k$, it meets a succession of $AC_i$-arcs and $BC_i$-arcs. A {\em trace} of ${\cal C}$ is the sequence obtained from $X$ by recording $0$ for each $AC_i$-arc included in $A' [x_0, y_0]$, $1$ for each $BC_i$-arc included in $P$ and $2$ for each $AC_i$-arc included in $A' [u_1, u_2]$. (This sequence is not uniquely determined unless $A$ and $B$ have been specified.) The next two claims show that, up to homeomorphism, the graph $$G [A \cup B \cup \bigcup_{i=1}^k C_i]$$ is determined by a trace of ${\cal C}$.

\medskip {\bf Claim 2:} Suppose that for some $i$ and $j > i$ there are a $BC_i$-arc included in $P$ and a $BC_j$-arc included in $P$. Let the former arc join vertices $u$ and $v$, where $u \in VC'_i [x_i, v]$, and the latter arc vertices $w$ and $x$, where $w \in VC'_j [x_j, x]$. Then $x \in VP [w, y_0]$, $u \in VP [x, y_0]$ and $v \in VP [u, y_0]$.

By Lemma 4 we see that $v \in VP [u, y_0]$ and $x \in VP [w, y_0]$. We must show that $u \in VP [x, y_0]$. Suppose therefore that $w \in VP [v, y_0]$. We may suppose $u$, $v$, $w$, $x$ to have been chosen to minimise $|P [v, w]|$.  It follows that the circuit $$C'_i [x_i, v] \cup P [v, w] \cup C'_j [w, y_j] \cup A' [y_j, x_i]$$ contradicts the choice of $C_i$. The claim has now been verified. \medskip

The same argument can also be applied to $AC_i$-arcs and $AC_j$-arcs meeting $A' [u_1, u_2]$. In other words, suppose that for some $i$ and $j > i$ there are an $AC_i$-arc meeting $A' [u_1, u_2]$ and an $AC_j$-arc meeting $A' [u_1, u_2]$. Let the former arc join vertices $u$ and $v$, where $u \in VC'_i [x_i, v]$, and the latter arc vertices $w$ and $x$, where $w \in VC'_j [x_j, x]$. Then $x \in VA' [w, u_2]$, $u \in VA' [x, u_2]$ and $v \in VA' [u, u_2]$.

\medskip {\bf Claim 3:} Suppose that for some $i$ there are two $BC_i$-arcs included in $P$. Let one of these arcs join vertices $u$ and $v$, where $u \in VC'_i [x_i, v]$, and let the other join vertices $w$ and $x$, where $w \in VC'_i [v, y_i]$ and $x \in VC'_i [w, y_i]$. Then $x \in VP [w, y_0]$, $u \in VP [x, y_0]$ and $v \in VP [u, y_0]$. 

By Lemma 4 we have $v \in VP [u, y_0]$ and $x \in VP [w, y_0]$. It remains to show that $u \in VP [x, y_0]$. Suppose therefore that $w \in VP [v, y_0]$. We may suppose $u$, $v$, $w$, $x$ to have been chosen to minimise $|P [v, w]|$. It follows from Claim 2 that we may now construct from ${\cal C}$ a new cascade ${\cal C'}$ by replacing $C_i$ with the circuit $$C'_i [x_i, v] \cup P [v, w] \cup C'_i [w, y_i] \cup A' [y_i, x_i].$$ This result contradicts the minimality of (1), and thereby completes the proof of the claim. \medskip

The same argument can be applied to $AC_i$-arcs which meet $A' [u_1, u_2]$. In other words, let us suppose that for some $i$ there are two such $AC_i$-arcs. Let one of these arcs join vertices $u$ and $v$, where $u \in VC'_i [x_i, v]$, and let the other join vertices $w$ and $x$, where $w \in VC'_i [v, y_i]$ and $x \in VC'_i [w, y_i]$. Then $x \in VA' [w, u_2]$, $u \in VA' [x, u_2]$ and $v \in VA' [u, u_2]$. 

We now investigate the possible traces a cascade may have. 

\medskip {\bf Claim 4:} Any two consecutive digits in a trace of ${\cal C}$ are distinct.

Suppose that a trace of ${\cal C}$ contains consecutive ones. Then for some $i$ there exist a $BC_i$-arc included in $P$ and joining vertices $u$ and $v$, where $u \in VC'_i [x_i, v]$, and another $BC_i$-arc included in $P$ and joining vertices $w$ and $x$, where $w \in VC'_i [v, y_i]$ and $x \in VC'_i [w, y_i]$, such that $C'_i [v, w]$ has no edges or internal vertices in common with the graph $G [A \cup B]$. By Claim 2 we have $x \in VP [w, y_0]$, $u \in VP [x, y_0]$ and $v \in VP [u, y_0]$. Moreover the edges of $f$ incident on $v$ and on $w$ both belong to $P [v, w]$. Let
\begin{eqnarray}
C &=& P [v, w] \cup C'_i [w, v], \nonumber \\
f^* &=& f+C, \nonumber \\
B^* &=& B+C \nonumber \\
&=& A' [y_0, u_1] \cup \{e_1\} \cup A' [v_1, x_0] \cup P [x_0, w] \cup C'_i [w, v] \cup P [v, y_0] \nonumber
\end{eqnarray}
and
\begin{eqnarray}
C^*_i &=& C_i + C \nonumber \\
&=& C'_i [x_i, u] \cup P [u, x] \cup C'_i [x, y_i] \cup A' [y_i, x_i]. \nonumber 
\end{eqnarray}
We now obtain a new cascade from ${\cal C}$ by replacing $C_i$ with $C^*_i$. This result contradicts the minimality of (1). 

The argument is similar if a trace of ${\cal C}$ contains consecutive twos. Suppose therefore that it contains consecutive zeros. Then for some $i$ such that $1 < i < k$ the circuit $C_i$ has a unique $AC_i$-arc and no $BC_i$-arc. The $AC_i$-arc joins vertices $x_i$ and $y_i$. Let 
\begin{eqnarray}
f^* &=& f+C_i, \nonumber \\
A^* &=& A+C_i \nonumber \\
&=& A' [v_1, y_i] \cup C'_i [y_i, x_i] \cup A' [x_i, u_1] \cup \{e_1\} \nonumber
\end{eqnarray}
and
\begin{eqnarray}
C^*_i &=& C_{i-1} + C_i + C_{i+1} \nonumber \\
&=& C'_{i-1} [x_{i-1}, y_{i-1}] \cup A' [y_{i-1}, x_{i+1}] \cup C'_{i+1} [x_{i+1}, y_{i+1}] \nonumber \\
& & {}\cup A' [y_{i+1}, y_i] \cup C'_i [y_i, x_i] \cup A' [x_i, x_{i-1}]. \nonumber
\end{eqnarray} 
Then $f^*$, $A^*$, $B$ and the cascade $$(C_1, C_2,\cdots, C_{i-2}, C^*_i, C_{i+2}, C_{i+3},\cdots, C_k)$$ contradict the minimality of (1). The proof of the claim is now complete.

\medskip {\bf Claim 5:} A trace of ${\cal C}$ contains neither $01$ nor $10$.

We may assume without loss of generality that a trace of ${\cal C}$ contains $01$. Then there exists $i>1$ such that $C'_i$ includes a path joining $x_i$ to a vertex $x \in VP - \{x_0, y_0\}$ and having no edges or internal vertices in common with the graph $G [A \cup B]$. By Lemma 4 the edge of $f$ incident on $x$ belongs to $P [x, y_0]$. Observe that the circuit $$Y = C'_i [x_i, x] \cup P [x, y_0] \cup A' [y_0, x_i]$$ is clockwise odd under a Pfaffian orientation of $K$, for otherwise the choice of $A$ and $B$ is contradicted by $A$ and $$B^* = A' [v_1, x_i] \cup C'_i [x_i, x] \cup P [x, y_0] \cup A' [y_0, u_1] \cup \{e_1\}$$ since $A' [x_0, x_i] \subseteq A \cap B^*$. Similarly the circuit $$Z = C'_i [x_i, x] \cup P [x, x_0] \cup A' [x_0, x_i]$$ is clockwise odd, as otherwise the choice of $A$ and $B$ is contradicted by $B^*$ and $B$ since $P[x, y_0] \subseteq B \cap B^*$. Thus of the three circuits $Y$, $Z$ and $$Y+Z = A+B = P \cup A' [x_0, y_0],$$ which have empty sum, only $A+B$ is clockwise even, yet an even number of them are clockwise even under a given orientation of $G$. This contradiction proves the claim. \medskip

By Claims 4 and 5 every second member of a trace of ${\cal C}$ must be 2. Having established the existence of a cascade with this property, we may now drop the requirement that $A$ and $B$ be chosen to maximise $|A' [v_2, y_0]| + |A' [x_0, v_1]|$. We obtain a {\em reduced} trace of ${\cal C}$ by removing each 2 from a trace. We deal next with the case where a reduced trace has three consecutive equal digits.

\medskip {\bf Claim 6:} Suppose a trace of ${\cal C}$ contains $02020$ or $12121$. Then $G$ has an even subdivision $H$ of $K_{3, 3}$ such that $G-VH$ has a $1$-factor. 

Without loss of generality we may assume that a trace of ${\cal C}$ contains $02020$. Then for some $i$ such that $2 \le i \le k - 2$ the circuit $C_i$ has a unique $AC_i$-arc meeting $A' [u_1, u_2]$ but it has no $BC_i$-arc meeting $P$, and a similar statement holds for $C_{i+1}$. The required subdivision of $K_{3, 3}$ is $$G [A' [u_1, y_0] \cup A' [x_0, v_1] \cup \{e_1\} \cup P \cup C_i \cup C_{i+1}].$$ 
\newline

We may henceforth assume that a reduced trace does not have three consecutive equal digits. 

A set $S$ of digits in a reduced trace is said to be {\em separated} if no two elements of $S$ are consecutive. We now consider the case where a reduced trace has three separated equal digits. We shall use $*$ to denote a non-empty subsequence of a trace of ${\cal C}$.

\medskip {\bf Claim 7:} Suppose a trace of ${\cal C}$ contains $12*212*21$ or $02*202*20$. Then $G$ has an even subdivision $H$ of $K_{3, 3}$ such that $G-VH$ has a $1$-factor. 

Without loss of generality we may suppose that a trace of ${\cal C}$ contains $02*202*20$. For some $h$ such that $2 \le h \le k - 2$ the path $C'_h$ includes a subpath joining $x_h$ to a vertex $x'_h \in VA' [u_1, u_2]$ and having no edges or internal vertices in common with the graph $G [A \cup B]$. Moreover for some $i$ such that $h < i < k$ the path $C'_{i-1}$ includes a subpath joining $y_{i-1}$ to a vertex $y'_{i-1} \in VA' [u_1, x'_h]$ and having no edges or internal vertices in common with $G [A \cup B]$, and $C'_i$ includes a subpath joining $x_i$ to a vertex $x'_i \in VA' [u_1, y'_{i-1}]$ and having no edges or internal vertices in common with $G [A \cup B]$. Furthermore, for some $j$ such that $i \le j < k$ the path $C'_j$ includes a subpath joining $y_j$ to a vertex $y'_j \in VA' [u_1, x'_i]$ and having no edges or internal vertices in common with $G [A \cup B]$. The required subdivision of $K_{3, 3}$ is 

\begin{center}
$G [A' [u_1, y_0] \cup A' [x_h, y_j] \cup A' [x_0, v_1] \cup \{e_1\} \cup P \cup C'_h [x_h, x'_h]$ \\
$\cup ~  C'_{i-1} [y_{i-1}, y'_{i-1}] \cup C'_i [x_i, x'_i] \cup C'_j [y_j, y'_j]].$
\end{center}
 
We may henceforth assume that a reduced trace does not have three separated equal digits.

At this juncture we observe that there are only a finite number of possible traces left to investigate, because any reduced trace of length at least $9$ must contain three consecutive zeros, three consecutive ones, three separated zeros or three separated ones. It can be shown that in each case there is a set of $f$-alternating circuits with sum $A + B$. If the clockwise parities of these circuits under a Pfaffian orientation of $K$ do not yield a contradiction, then the fact that $A + B$ is clockwise even under a Pfaffian orientation of $K$ implies, by Lemma 2, that exactly one of $A$ and $B$ is clockwise even. In this case $G$ must be non-Pfaffian, as it cannot be oriented so that all the $f$-alternating circuits are clockwise odd. In such graphs $G$ we therefore search for a subgraph $H$, reducible to an even subdivision of $K_{3, 3}$, such that $G-VH$ has a $1$-factor.

The next claim is proved by investigating the condition under which a contradiction can be derived from the clockwise parities, under a Pfaffian orientation of $K$, of a set of $f$-alternating circuits with sum $A+B$.

\medskip {\bf Claim 8:} A reduced trace of ${\cal C}$ is of odd length.

Suppose that a reduced trace of ${\cal C}$ is of even length. First we construct a cascade ${\cal D}$ along $P$ from $y_0$ to $x_0$. If a reduced trace has no $1$, then ${\cal D}$ consists only of the circuit $$D_1 = X \cup A' [y_k, x_0] \cup P \cup A' [y_0, x_1].$$ Otherwise there exists a vertex $z_1 \in VX \cap VP$. We may assume $z_1$ chosen to minimise $|X [w_1, z_1]|$, where $w_1 = x_1$. In this case we define $$D_1 = X[w_1, z_1] \cup P [z_1, y_0] \cup A' [y_0, w_1].$$ We also let $w_2$ be the vertex of $VX \cap (VP [z_1, z_0] - \{z_0\})$ that minimises $|P [w_2, z_0]|$, where $z_0 = y_0$. 

If there exists a vertex $z_2 \in (VX [w_2, y_k] - \{w_2\}) \cap VP$, then we repeat the argument. Thus we may assume $z_2$ chosen to minimise $|X [w_2, z_2]|$. We then define $$D_2 = X [w_2, z_2] \cup P [z_2, w_2].$$ We also let $w_3$ be the vertex of $VX \cap (VP [z_2, z_1] - \{z_1\})$ that minimises $|P [w_3, z_1]|$.

If there exists a vertex $z_3 \in (VX [w_3, y_k] - \{w_3\}) \cap VP$, then we repeat the argument inductively. Thus there exists $l > 0$ such that $$(VX [w_l, y_k] - \{w_l\}) \cap VP = \emptyset.$$ We define $$D_l = X [w_l, y_k] \cup A' [y_k, x_0] \cup P [x_0, w_l].$$ The required cascade ${\cal D}$ is then $(D_1, D_2,\cdots, D_l)$. 

We observe that the number of ones in a trace of ${\cal C}$ is $l - 1$. Similarly the number of zeros is $k - 1$. Hence $k + l$ is even by assumption. Moreover $$\sum^k_{i = 1} C_i + \sum^l_{j = 1} D_j = A + B.$$ (For any $i$ each member of a $PD_i$-arc belongs to one circuit in ${\cal C}$ and two circuits in ${\cal D}$. The remaining elements of $P$ belong to just one circuit in ${\cal D}$ and to no circuits in ${\cal C}$. Similar results hold for the edges of $A' [x_0, y_0]$. The remaining edges of $X \cup A' [y_0, x_1] \cup A' [y_k, x_0]$ belong to one circuit in each of ${\cal C}$ and ${\cal D}$.) In addition, it is easy to see that $K$ may be oriented so that the circuits in ${\cal C}$ and ${\cal D}$ are directed. For example, if $\{M, N\}$ is a bipartition of $K$, we may orient the edges of $f$ from $M$ to $N$ and the remaining edges from $N$ to $M$. Under this orientation the circuit $A + B$ is clockwise odd. As it is the only clockwise odd circuit in the set $$S = \{C_1, C_2,\cdots, C_k, D_1, D_2,\cdots, D_l, A+B\},$$ and $S$ has empty sum and odd cardinality, we now have a contradiction because under a Pfaffian orientation of $K$ the circuit $A + B$ is the only clockwise even member of $S$. This contradiction proves the claim. \medskip

In all the remaining cases we must find a subgraph $H$ of $G$, reducible to an even subdivision of $K_{3, 3}$, such that $G - VH$ has a $1$-factor. We simplify the work by establishing two further claims which enable us to manipulate a trace of ${\cal C}$.

\medskip {\bf Claim 9:} Let $G'$ be a graph with a subgraph $H'$, reducible to an even subdivision of $K_{3, 3}$, such that $G' - VH'$ has a $1$-factor. Suppose that $G'$ has a cascade ${\cal C'}$ with a trace whose first element is $0$ or $1$. Let $G$ have a cascade ${\cal C}$ with a trace consisting of $2$ followed by the trace of ${\cal C'}$. Then $G$ has a subgraph $H$, reducible to an even subdivision of $K_{3, 3}$, such that $G-VH$ has a $1$-factor.

The path $C'_1$ includes a subpath joining $x_1$ to a vertex $x'_1 \in VA' [u_1, u_2]$ and having no edges or internal vertices in common with $G[A \cup B]$. Contraction of the circuit $C'_1 [x_1, x'_1] \cup A' [x'_1, x_1]$ yields the graph $G'$. Since $G'$ has a subgraph $H'$, reducible to an even subdivision of $K_{3, 3}$, such that $G'-VH'$ has a $1$-factor, the claim follows. \medskip

A corresponding result holds if the last element of a trace of ${\cal C'}$ is $0$ or $1$ and a trace of ${\cal C}$ consists of that of ${\cal C'}$ followed by $2$. We may therefore assume that neither the first nor the last element of a trace of ${\cal C}$ is $2$.

\medskip {\bf Claim 10:} Suppose that the first element of a trace of ${\cal C}$ is $0$. Then ${\cal C}$ has another trace obtained by replacing the first $0$ with $1$.

Note first that the path $C'_1 [x_1, y_1]$ has no edges or internal vertices in common with $G [A \cup B]$. Let 
\begin{eqnarray}
f^* &=& f + C_1, \nonumber \\
A^* &=& A + C_1 \nonumber \\
&=& A' [u_1, x_1] \cup C'_1 [x_1, y_1] \cup A' [y_1, v_1] \cup \{e_1\} \nonumber
\end{eqnarray}
and
\begin{eqnarray}
B^* &=& B + C_1 \nonumber \\
&=& P \cup A' [y_0, y_1] \cup C'_1 [y_1, x_1] \cup A' [x_1, u_1] \cup \{e_1\} \cup A' [v_1, x_0]. \nonumber
\end{eqnarray}
The required trace is calculated by replacing $A$ and $B$ with the $f^*$-alternating circuits $A^*$ and $B^*$ respectively. \medskip

A corresponding result holds if the first element of a trace of ${\cal C}$ is $1$. We also obtain corresponding results for the last element of a trace of ${\cal C}$.

\medskip {\bf Claim 11:} A reduced trace of ${\cal C}$ has length $1$.

Suppose first that a reduced trace of ${\cal C}$ has length at least $5$. By Claim 10 we may assume that the first and last digits are $0$. By Claim $7$ the third digit must therefore be $1$. However, by Claim 10 there is another reduced trace of ${\cal C}$ obtained by replacing the first and last digits with $1$, in violation of Claim 7. The argument is similar if a reduced trace of ${\cal C}$ has length $3$, except that Claim 6 is used instead of Claim 7. Claim 11 now follows from Claim 8. \medskip

Thus ${\cal C}$ has a reduced trace of length $1$. By symmetry and Claim 9, we may assume that $0$ is a trace of ${\cal C}$. Then $$G [A \cup B \cup C_1 \cup C_2]$$ is the required subdivision of $K_{3, 3}$. \epr
\newline
\newline
{\bf Theorem 5}  {\it Let $G$ be a $1$-extendible graph that
  can be obtained from a $1$-extendible bipartite graph $K$ by a
  $2$-ear adjunction. Suppose that $G [C \cup D]$ is bipartite
  whenever $C$ and $D$ are alternating circuits that traverse both
  ears and have the property that there are just two $\bar CD$-arcs.
  Then $G$ is non-Pfaffian if and only if it has a subgraph $H$,
  reducible to an even subdivision of $K_{3,3}$, such that $G-VH$ has
  a $1$-factor.}  \newline

{\bf Proof.} We have already seen that if $G$ has $H$ as a subgraph
then $G$ is non-Pfaffian. Suppose therefore that $G$ is non-Pfaffian.
Suppose also that $G$ is formed from a $1$-extendible bipartite graph
$K$ by the adjunction of ears $E_1$ and $E_2$. If $K$ is non-Pfaffian,
then $K$ has the required subgraph $H$ by Theorem $2$, and therefore
so does $G$. We may suppose therefore that $K$ is Pfaffian.

Suppose that $G-E_1$ is $1$-extendible. Since $G-E_1$ is obtained from
the bipartite graph $K$ by the adjunction of the single ear $E_2$, it
follows that $G-E_1$ is also bipartite. Since $G$ is obtained from
$G-E_1$ by the adjunction of a single ear, $G$ is also bipartite.
Therefore $G$ has the required subgraph by Theorem $2$. The argument
is similar if $G-E_2$ is $1$-extendible.

We may now assume that neither $G-E_1$ nor $G-E_2$ is $1$-extendible.
Let $E_1$ join vertices $u_1$ and $v_1$ and fix a $1$-factor $f$ of
$K$. Since $K$ is $1$-extendible, it has a path $P$, joining $u_1$ and
$v_1$, such that each vertex of $VP-\{u_1\}$ is incident with an edge
of $P \cap f$. Thus $P \cup E_1$ is an alternating circuit
in $G-E_2$ if $G-E_2$ is bipartite, and in this case we reach the
contradiction that $G-E_2$ is $1$-extendible. Hence $G-E_2$ is not
bipartite. Similarly $G-E_1$ is not bipartite.

Next, we may clearly assume without loss of generality that $E_1$ and
$E_2$ both have cardinality $1$. Therefore we may let $E_1=\{e_1\}$
and $E_2=\{e_2\}$. Thus $e_1$ and $e_2$ are distinct edges of $EG-f$
such that neither $G-\{e_1\}$ nor $G-\{e_2\}$ is bipartite and neither
is $1$-extendible, but $G-\{e_1, e_2\}$ is bipartite and
$1$-extendible. As $G$ is $1$-extendible, there must be an
$f$-alternating circuit $X$ containing $e_1$. Since $G-\{e_2\}$ is not
$1$-extendible, $X$ must also contain $e_2$. It follows that $e_1$ and
$e_2$ are independent edges, as $X$ is $f$-alternating.

Since $K$ is Pfaffian we may assume given a Pfaffian orientation of
$K$. Extend this orientation to $G$ by orienting $e_1$ and $e_2$
arbitrarily. As $G$ is not Pfaffian, there must exist a clockwise even
alternating circuit $A$ in $G$. Recalling that some $1$-factor in $G$
must have its sign opposite to that of $f$, we find that $A$ may be
chosen to be $f$-alternating. Being clockwise even, it cannot be
alternating in $K$, and therefore cannot be a circuit of $K$. Thus $A$
contains $e_1$ or $e_2$ and hence both. There must also be a clockwise
odd alternating circuit $B$ containing $e_1$ and $e_2$ for otherwise a
Pfaffian orientation of $G$ could be realised by reorienting $e_1$ or
$e_2$. We can also choose $B$ to be $f$-alternating. By Lemma $3$ we
find that $f$, $A$ and $B$ may be chosen so that there are just one or
two $AB$-arcs, each containing $e_1$ or $e_2$ and their union
containing both.

If there are two $AB$-arcs, then $G[A \cup B]$ is bipartite by
hypothesis. However one of these arcs must contain $e_1$ and the other
$e_2$. Let $A'=A-\{e_1\}$, and let $e_2$ join vertices $u_2$ and
$v_2$, where $e_2 \in A' [u_1, v_2]$. If there exists an $\bar
AB$-arc $Q$ joining a vertex $x \in VA' [u_1, u_2]$ to a vertex $y
\in VA' [v_1, v_2],$ then since $G [A \cup B]$ is bipartite we find
that either $$Q \cup A' [y, x]$$
or $$Q \cup A' [x, u_1] \cup \{e_1\}
\cup A' [v_1, y]$$
is an $f$-alternating circuit containing just one
member of $\{e_1, e_2\}$. From this contradiction to the fact that
neither $G-\{e_1\}$ nor $G-\{e_2\}$ is $1$-extendible, together with
the knowledge that one $AB$-arc contains $e_1$ and the other $e_2$, we
deduce that $A+B$ is the union of two disjoint circuits.

By Lemma $2$ one of these circuits is a clockwise even alternating
circuit. It is therefore possible to modify $B$ so that there is only
one $AB$-arc. Then $A+B$ is a clockwise even alternating circuit. We can now invoke Theorem $4$ to draw the
desired conclusion. \epr \newline

In order to complete the analysis of a $2$-ear adjunction to a bipartite graph, it remains to consider the case where there are two $AB$-arcs. This case is considerably more complicated, and is currently being studied by the authors.

\end{document}